\documentclass[a4paper]{article}
\usepackage{amsfonts}
\usepackage{amsthm}
\usepackage{amsmath}
\numberwithin {equation}{section}
\newcommand{\N}{\mathbb{N}\sb{0}}
\newcommand{\R}{\mathbb{R}\sp{n}}

\newtheorem{theorem}{Theorem}[section]
\newtheorem{lemma}[theorem]{Lemma}

\newtheorem{corollary}[theorem]{Corollary}

\newtheorem{definition}[theorem]{Definition}

\newtheorem{remark}[theorem]{Remark}

\newcommand{\q}{q(x,D)u(x):=(2\pi)^{-\frac{n}{2}}\int_{\R} e^{ix\cdot \xi} q(x,\xi)\hat u(\xi) d\xi}

\newcommand{\n}{\newline}
\newcommand{\ps}{pseudo-differential operator}
\newcommand{\parxi}{\partial^{\alpha}_{\xi}}
\newcommand{\parx}{\partial^{\beta}_{x}}
\newcommand{\Sr}{S^{m,\psi}_{\rho}(\R)}
\newcommand{\Sz}{S^{m,\psi}_{0}(\R)}

\begin{document} 
{\begin{center}
{\Huge {Feller Semigroups Obtained by Variable Order Subordination}}
\end{center}}
\vskip15pt
{\begin{center}
{\Large{Kristian P. Evans and Niels Jacob}}
\end{center}}
\begin{center}
{\Large{July 2006}}
\end{center}
\section*{Abstract}
For certain classes of negative definite symbols $q(x,\xi)$ and state space dependent Bernstein function $f(x,s)$ we prove that $-p(x,D)$, the pseudo-differential operator with symbol $-p(x,\xi)=-f(x,q(x,\xi))$, extends to the generator of a Feller semigroup. Our result extends previously known results related to operators of variable (fractional) order of differentiation, or variable order fractional powers. New concrete examples are given.
\section{Introduction}
In the early days of the theory of pseudo-differential operators, pseudo differential operators of variable order had already been studied, compare A. Unterberger and J. Bokobza [\ref{AAV}]. These considerations were taken up by H.-G. Leopold [\ref{AAQ}], [\ref{AAR}] who gave more emphasis on the function space point of view. On the other hand, also in the early days of the theory of pseudo-differential operators Ph. Courr\`{e}ge [\ref{AAB}] pointed out that (most) generators of Feller semigroups are pseudo-differential operators, but their symbols do not belong to``nice'' or``classical'' symbol classes. Indeed, on $S(\R)$ the generator of a Feller semigroup has the representation
\begin{equation}
Au(x) = -q(x,D)u(x) = -(2\pi)^{-\frac{n}{2}}\int_{\R} e^{ix\cdot \xi} q(x,\xi)\hat u(\xi) d\xi
\end{equation}
where the symbol $q:\R \times \R \rightarrow \mathbb{C}$ is measurable and locally bounded and for $x \in \R$ fixed $q(x, \cdot)$ is a continuous negative definite function, i.e. we have the L\`{e}vy-Khinchin representation
\begin{equation}
q(x,\xi)=c(x) +id(x)\xi+\sum_{k,l=1}^na_{k,l}(x)\xi_k\xi_l+\int_{\R\setminus\{0\}} \left(1-e^{-iy\cdot \xi}-\frac{iy\cdot \xi}{1+|y|^2}\right)\nu(x,dy)
\end{equation}
with $c(x)\geq0$, $d(x)\in \R$, $a_{kl}(x)=a_{lk}(x) \in \mathbb{R}$ and $\sum_{k,l=1}^na_{kl}(x)\xi_k\xi_l\geq 0$, and $\int_{\R\setminus\{0\}}(1 \wedge |y|^2)\nu(x,dy) < \infty$. Thus these symbols need not to be smooth with respect to $\xi$ nor do they need to have a nice expansion into homogeneous functions. Maybe the fact that these symbols are a bit exotic is the reason why Courr\`{e}ge's result was almost ignored for around 25 years. In [\ref{AAK}], see also [\ref{AAJ}], Courr\`{e}ge's idea was taken up and a systematic study of pseudo-differential operators generating Markov processes was initiated, see also [\ref{AAL}] - [\ref{AAN}]. 
\par The fact that the composition of a Bernstein function $f$ with a continuous negative definite function $\psi$ is again a continuous negative definite function gives a powerful tool to construct new (Feller) semigroups from given ones. If $q(x,\xi)$ is a suitable symbol such that $-q(x,D)$ generates a Feller semigroup, then $(f\circ q)(x,\xi)=f(q(x.\xi))$ is a symbol with the property that $\xi \rightarrow (f\circ q)(x,\xi)$ is a continuous negative definite function and therefore $-(f\circ q)(x,D)$ is a candidate for being a generator of a Feller semigroup. Of course, this procedure is closely linked to subordination in the sense of Bochner. 
\par In a joint paper [\ref{AAO}] with H. G. Leopold it was suggested to study Feller semigroups obtained by subordination of variable order, more precisely, to consider ``fractional powers of variable order'' in case of the symbol $(1+|\xi|^2)$, i.e. to study $(x,\xi) \rightarrow (1+|\xi|^2)^{\alpha(x)}$. These ideas were taken up and further investigations on fractional powers of variable order are due to A. Negoro [\ref{AAU}], K. Kikuchi and A. Negoro [\ref{AAP}], as well as F. Baldus [\ref{AAA}]. Finally, W. Hoh in [\ref{AAG}] could combine his symbolic calculus [\ref{AAE}] with these ideas, compare W. Hoh [\ref{AAF}] and [\ref{AAH}].
\par The purpose of this note is twofold. First we suggest a method to study ``variable order subordination'' for more general Bernstein functions than \newline $f_{\alpha}(s)=s^{\alpha}$, $0 < \alpha < 1$. More precisely, we consider symbols of the form 
\begin{equation}
p(x,\xi)=f(x,q(x,\xi))
\end{equation}
where $q$ is a suitable symbol from Hoh's class and $f:\R\times [0,\infty) \rightarrow \mathbb{R}$ is a smooth function such that for fixed $x \in \R$ the function $s\rightarrow f(x,s)$ is a Bernstein function. Our method uses some ideas from the theory of t-coercive (differential) operators as investigated by I. S. Louhivaara and C. Simader [\ref{AAS}]-[\ref{AAT}] in order to establish the result that $-p(x,D)$ generates a Feller semigroup. Secondly, we enrich the class of examples by studying the Bernstein function 
$$
s \rightarrow s^{\frac{\alpha}{2}}(1 - e^{-4s^{\frac{\alpha}{2}}}).
$$
Since we depend on Hoh's symbolic calculus we recollect some basic facts of this calculus in our first section. All our methods are standard, i.e. they are as in [\ref{AAL}]-[\ref{AAN}].
\vskip15pt
The first named author is grateful for support obtained from Swansea University and EPSRC.
\section{Hoh's Symbolic Calculus}
Before starting with our main considerations we need to recollect some basic results from Hoh's symbolic calculus, see W.Hoh [\ref{AAE}] or [\ref{AAF}], compare also [\ref{AAM}].
\begin{definition}
A continuous negative definite function $\psi:\R \rightarrow \R$ belongs to the class $\Lambda$ if for all $\alpha \in \N^n$ it satisfies
\begin{equation}
|\parxi(1+\psi(\xi))| \leq c_{|\alpha|}(1+\psi(\xi))^{\frac{2-\rho(|\alpha|)}{2}},
\end{equation}
where $\rho(k)=k \wedge 2$ for $k \in \N^n$.
\end{definition}
\begin{definition}
A. Let $m \in \mathbb{R}$ and $\psi \in \Lambda$. We then call a $C^{\infty}$-function $q:\R \times \R \longrightarrow \mathbb{C}$ a symbol in the class $S^{m,\psi}_{\rho}(\R)$ if for all $\alpha,\beta \in \N^n$ there are constants $c_{\alpha,\beta}\geq 0$ such that
\begin{equation}
|\parx \parxi q(x,\xi)| \leq c_{\alpha,\beta}(1+\psi(\xi))^{\frac{m-\rho(|\alpha|)}{2}}
\end{equation}
holds for all $x\in \R$ and $\xi \in \R$. We call $m \in \mathbb{R}$ the order of the symbol $q(x, \xi)$.
\n 
B. Let $\psi \in \Lambda$ and suppose that for an arbitrarily often differentiable function $q:\R \times \R \longrightarrow \mathbb{C}$ the estimate
\begin{equation}
|\parxi \parx q(x,\xi)| \leq \tilde{c}_{\alpha,\beta}(1+\psi(\xi))^{\frac{m}{2}}
\end{equation}
holds for all $\alpha, \beta \in \N^n$ and $x, \xi \in\R$. In this case we call $q$ a symbol of the class $S^{m,\psi}_{0}(\R)$.
\end{definition}
Note that $S^{m,\psi}_{\rho}(\R) \subset S^{m,\psi}_{0}(\R)$.
For $q\in \Sz$, hence also for $q\in \Sr$,  we can define on $S(\R)$ the \ps\ $q(x,D)$ by 
\begin{equation}
\q
\end{equation}
and we denote the classes of these operators by $\Psi^{m,\psi}_{\rho}(\R)$ and $\Psi^{m,\psi}_{0}(\R)$, respectively.
\begin{theorem}
Let $q\in \Sz$ then $q(x,D)$ maps $S(\R)$ continuously into itself.
\end{theorem}
Let $\psi:\R \rightarrow \mathbb{R}$ be a fixed continuous negative definite function. For $s \in \mathbb{R}$ and $u \in S(\R)$ (or $u \in S'(\R)$) we define the norm 
\begin{equation}
||u||^2_{\psi,s} = ||(1+\psi(D))^{\frac{1}{2}}u||^2_{0}=\int_{\R}(1+\psi(s))^s|\hat{u}(\xi)|^2d\xi.
\end{equation}
The space $H^{\psi,s}(\R)$ is defined as 
\begin{equation}
H^{\psi,s}(\R):=\{u \in S'(\R); ||u||_{\psi,s} <\infty\}.
\end{equation}
The scale $H^{\psi,s}(\R)$, $s \in \R$, and more general spaces have been systematically investigated in [\ref{AAC}] and [\ref{AAD}], see also [\ref{AAM}]. In particular we know that if for some $\rho_1>0$ and $\tilde{c}_1 >0$ the estimate $\psi(\xi) \geq \tilde{c}_1|\xi|^{\rho_1}$ holds for all $\xi \in \R$, $|\xi| \geq R, \ R\geq 0$, then the space $H^{\psi,s}(\R)$ is continuously embedded into $C^{\infty}(\R)$ provided $s> \frac{n}{2 \rho_1}$.
\begin{theorem}
Let $q\in \Sz$ and let $q(x,D)$ be the corresponding \ps. For all $s \in \mathbb{R}$ the operator $q(x,D)$ maps the space \newline$H^{\psi,m+s}(\R)$ continuously into the space $H^{\psi,s}(\R)$, and for all $u \in H^{\psi, m+s}(\R)$ we have the estimate
\begin{equation}
||q(x,D)u||_{\psi,s} \leq c||u||_{\psi,m+s}.
\end{equation}
\end{theorem}
On $S(\R)$ we may define the bilinear form
\begin{equation}
B(u,v):=(q(x,D)u,v)_0,\ \ \ q \in \Sr.
\end{equation}
\begin{theorem}
Let $q \in \Sr$ be real valued and $m>0$. It follows that
\begin{equation}
|B(u,v)| \leq c||u||_{\psi, \frac{m}{2}}||v||_{\psi,\frac{m}{2}}
\end{equation}
holds for all $u,v \in S(\R)$. Hence the bilinear form $B$ has a continuous extension onto $H^{\psi,\frac{m}{2}}(\R)$. If in addition for all $x \in \R$
\begin{equation}
q(x,\xi)\geq \delta_0(1+\psi(\xi))^{\frac{m}{2}} \mbox{ for } |\xi| \geq R
\label{AA}
\end{equation}
with some $\delta_0>0$ and $R \geq 0$, and 
\begin{equation}
\lim_{|\xi|\rightarrow \infty} \psi(\xi) =\infty
\label{AB}
\end{equation}
holds, then we have for all $u \in H^{\psi,\frac{m}{2}}(\R)$ the G{\aa}rding inequality
\begin{equation}
ReB(u,u) \geq \frac{\delta_0}{2}||u||^2_{\psi, \frac{m}{2}}-\lambda_0||u||^2_0.
\label{AC}
\end{equation}
\end{theorem}
Furthermore we have 

\begin{theorem}
If we assume (\ref{AA}) and (\ref{AB}) then for $s >-m$ we have 
\begin{equation}
\frac{\delta_0}{2}||u||_{\psi,m+s} \leq ||q(x,D)u||^2_{\psi,s}+||u||^2_{\psi,m+s-\frac{1}{2}}
\end{equation}
for $q \in \Sr$ real-valued and all $u \in H^{\psi,s+m}(\R)$.
\end{theorem}
From Theorem 1.5 and 1.6 one may deduce the following regularity result:
\begin{theorem}
Let $q \in \Sr$ be as in Theorem 1.6, $m\geq 1$. Further suppose that for $f \in H^{\psi,s}(\R)$, $s\geq0$, there exists $u \in H^{\psi, \frac{m}{2}}(\R)$ such that
\begin{equation}
B(u,\phi)=(f,\phi)_{L^2}
\end{equation}
holds for all $\phi \in H^{\psi,\frac{m}{2}}(\R)$ (or $\phi \in S(\R)$). Then $u$ belongs already to the space $H^{\psi, m+s}(\R)$.
\end{theorem}
So far we have used properties of symbols to establish mapping properties and estimates for operators. The real power of a symbolic calculus is that it reduces calculations for operators to calculations for symbols. The following result is most important for us
\begin{theorem}
Let $\psi \in \Lambda$. For $q_1 \in S^{m_1,\psi}_{\rho}(\R)$ and $q_2 \in S^{m_2,\psi}_{\rho}(\R)$ the symbol $q$ of the operator $q(x,D):=q_1(x,D) \circ q_2(x,D)$ is given by 
\begin{equation}
q(x,\xi)=q_1(x,\xi)\cdot q_2(x,\xi) + \sum_{j=1}^{n} \partial_{\xi_j}q_1(x,\xi) D_{x_j}q_2(x,\xi) +q_{r_1}(x,\xi)
\label{AD}
\end{equation}
with $q_{r_1} \in S^{m_1+m_2-2,\psi}_0(\R)$.
\end{theorem}
\begin{remark}
An easy calculation yields $q_1 \cdot q_2 \in S^{m_1+m_2,\psi}_{\rho}(\R)$, \n$\partial_{\xi_j}q_1 \in S^{m_1-1,\psi}_{\rho}(\R)$ and $D_{x_j}q_2 \in S^{m_2,\psi}_{\rho}(\R)$. Hence the second term on the right hand side in (\ref{AD}) belongs to $S^{m_1+m_2-1,\psi}_{\rho}(\R)$.
\end{remark}

\section{The Formal Background of our Proof that \n$-p(x,D)$ Generates a Feller Semigroup}
The proof that $-p(x,D)$ as described in the introduction, see also below, extends to a generator of a Feller semigroup depends on various estimates which might be different for different operators. However, once these estimates are established we only need to apply a piece of ``soft'' analysis. In this section we discuss this part of the proof, i.e. we will assume all crucial estimates hold.
\n Let $f: \R \times [0,\infty) \rightarrow \mathbb{R}$ be an arbitrarily often differentiable function such that for $y \in \R$ fixed the function $s\rightarrow f(y,s)$ is a Bernstein function. Moreover we assume

\begin{equation}
\inf _{y \in \R} f(y,s) \geq f_0(s) \mbox{ \ \ \  for all } s \in [0,\infty)
\label{LL}
\end{equation}
as well as
\begin{equation}
\sup_{y \in \R} f(y,s) \leq f_1(s) \mbox{ \ \ \  for all } s \in[0,\infty)
\label{M}
\end{equation}
where $f_0$ and $f_1$ are Bernstein functions.
For a given real-valued negative definite symbol $q(x,\xi)$ it follows that
\begin{equation}
p(y;x,\xi):=f(y,q(x,\xi))
\end{equation}
give rise to a further negative definite symbol by defining
\begin{equation}
p(x,\xi):=p(x;x,\xi).
\label{N}
\end{equation}
In case where $q(x,\xi)$ is comparable with a fixed continuous negative definite function $\psi$, i.e.
\begin{equation}
0<c_0 \leq \frac{q(x,\xi)}{\psi(\xi)} \leq c_1, \ \ \ \ c_1\geq1,
\label{P}
\end{equation}
for all $x \in \R$ and $\xi \in \R$, we find using Lemma 3.9.34.B in [\ref{AAL}]
\begin{equation}
p(x,\xi) \leq f(y_1,q(x,\xi)) \leq c_1 f_1(\psi(\xi))
\end{equation}
and we define
\begin{equation}
\psi_1(\xi):=c_1f_1(\psi(\xi)).
\label{AAB}
\end{equation}
Moreover it holds
$$
p(x,\xi)\geq f(y_0,q(x,\xi)) \geq \tilde{c}_0f_0(\psi(\xi))
$$
and we set 
\begin{equation}
\psi_0(\xi):= c_0'f_0(\psi(\xi)).
\label{AAC}
\end{equation}
Clearly, $\psi_0$ and $\psi_1$ are continuous negative definite functions. Later on we assume that for $|\xi|$ large
\begin{equation}
\psi(\xi) \geq \tilde{c}_1 |\xi|^{\rho_1}, \ \ \ \tilde{c}_1>0 \mbox{ \ and \ } \rho_1>0
\label{FF}
\end{equation}
holds as well as
\begin{equation}
f(y_0,s) \geq \tilde{c}_0 s^{\rho_0}, \ \ \tilde{c}_0>0 \mbox{ \ and \ } \rho_0>0.
\label{GG}
\end{equation}
This implies for $|\xi|$ large that
\begin{equation}
\psi_0(\xi) \geq \tilde{c}_2 |\xi|^{\rho_0 \rho_1}, \ \ \tilde{c}_2 >0,
\label{R}
\end{equation}
holds. Since $\psi_0(\xi) \leq \psi_1(\xi)$ we have 
\begin{equation}
H^{\psi_1,1}(\R) \hookrightarrow H^{\psi_0,1}(\R).
\end{equation}
We add the assumption that there exists $0 < \sigma < \frac{1}{2}$ such that 
\begin{equation}
(1+\psi_1)^{\frac{1}{2}}\in S^{1+\sigma,\psi_0}_{\rho}(\R).
\label{A}
\end{equation}
This will imply that
\begin{equation}
H^{\psi_0, m(1+\sigma)}(\R) \hookrightarrow H^{\psi_1,m}(\R)
\label{JJ}
\end{equation}
holds for $m\geq 0$. Further, (\ref{A}) implies that if $p_1(x,\xi)$ is any symbol belonging to $S^{m,\psi_1}_{\rho}(\R)$ then it also belongs to $S^{m(1+\sigma),\psi_0}_{\rho}(\R)$ which follows from

\begin{eqnarray*}
|\parxi \parx p_1(x,\xi)| & \leq & c_{\alpha,\beta} (1+\psi_1(\xi))^{\frac{m - \rho(|\alpha|)}{2}} \\ 
& \leq & \tilde{c}_{\alpha,\beta}(1+\psi_0(\xi))^{\frac{m - \rho(|\alpha|)(1+\sigma)}{2}} \\
& \leq & \tilde{c}_{\alpha,\beta}(1+\psi_0(\xi))^{\frac{(1+\sigma)m-\rho(|\alpha|)}{2}}.
\end{eqnarray*}

The \ps\ $q(x,D)$ has the symbol $q\in S_{\rho}^{2,\psi}(\R)$. We assume that the \ps\ $p(x,D)$, defined on $S(\R)$ by
$$
p(x,D)u(x)=(2\pi)^{-\frac{n}{2}}\int_{\R} e^{ix\cdot \xi} p(x,\xi)\hat u(\xi) d\xi 
$$
\begin{equation}
 =(2\pi)^{-\frac{n}{2}}\int_{\R} e^{ix\cdot \xi}f(x,q(x,\xi))\hat u(\xi)d\xi
\end{equation}
has a symbol $p \in S^{2+\tau_1, \psi_1}_{\rho}(\R)$ for some appropriate $\tau_1 \geq 0$.
This implies together with (\ref{A}) that the operator $p(x,D)$ is continuous from \n$H^{\psi_0,2+\tau_1+2\sigma+\tau_1 \sigma+s}(\R)$ to $H^{\psi_0,s}(\R)$, in particular it is continuous from \newline$H^{\psi_0,1}(\R)$ to $H^{\psi_0,-1-\tau_1-2\sigma-\tau_1\sigma}(\R)$.
\n With $p(x,D)$ we can associate the bilinear form
\begin{equation}
B(u,v):= (p(x,D)u,v)_0,\ \ u,v \in S(\R).
\end{equation}
Assuming the estimate
\begin{equation}
|B(u,v)| \leq \kappa||u||_{\psi_1,1}||v||_{\psi_1,1}, \ \ \kappa \geq0,
\end{equation}
to hold for all $u,v \in S(\R)$, we may extend $B$ to a continuous bilinear form on $H^{\psi_1,1}(\R)$. This extension is again denoted by $B$. For $u \in H^{\psi_1,1}(\R)$ we assume in addition
\begin{equation}
B(u,u) \geq \gamma||u||^2_{\psi_0,1} - \lambda_0||u||^2_0, \ \lambda_0 \geq 0, \gamma > 0.
\label{D}
\end{equation}
Following ideas from I.S. Louhivaara and Chr. Simader, [\ref{AAS}] and [\ref{AAT}], we consider an intermediate space associated with 
\begin{equation}
B_{\lambda_0}(u,v):=B(u,v)+\lambda_0(u,v)_0,
\end{equation}
namely the space $H^{p_{\lambda_0}}(\R)$ defined as a completion of $S(\R)$ (or $H^{\psi_1,1}(\R)$) with respect to the scalar product $B_{\lambda_0}$. Obviously we have 
\begin{equation}
H^{\psi_1,1}(\R)\hookrightarrow H^{p_{\lambda_0}}(\R) \hookrightarrow H^{\psi_0,1}(\R)
\label{B}
\end{equation}
in the sense of continuous embeddings. Moreover, by the Lax-Milgram theorem, for every $f\in \left( H^{p_{\lambda_0}}(\R)\right)^{*}$ exists a unique element $u \in H^{p_{\lambda_0}}(\R)$ satisfying
\begin{equation}
B_{\lambda_0}(u,v) = <f,v>
\label{C}
\end{equation}
for all $v \in H^{p_{\lambda_0}}(\R)$. This element we call the variational solution to the equation $p(x,D)u+\lambda_0u=f$.
\n From (\ref{B}) we derive
\begin{equation}
H^{\psi_0,-1}(\R)=\left(H^{\psi_0,1}(\R)\right)^{*} \hookrightarrow \left( H^{p_{\lambda_0}}(\R)\right)^{*},
\end{equation}
hence for $f \in H^{\psi_0,-1}(\R)$ there exists a unique $u \in H^{p_{\lambda_0}}(\R)$ satisfying (\ref{C}).  We claim now that for every $f \in H^{\psi_0,-1}(\R)$ there exists a unique $u \in H^{\psi_0,1}(\R)$ such that
 \begin{equation}
p_{\lambda_0}(x,D)u =p(x,D)u + \lambda_0u = f
\label{E}
\end{equation}
holds. Denote by $u \in H^{p_{\lambda_0}}(\R)$ the unique solution to (\ref{C}) for $f \in H^{\psi_0,-1}(\R)$ given and take a sequence $(u_k)_{k \in \mathbb{N}}$, $u_k \in S(\R)$, converging in  $H^{p_{\lambda_0}}(\R)$ to $u$. It follows from
$$
(p_{\lambda_0}(x,D)u_k,v)_0=B_{\lambda_0}(u_k,v), \ \ \ v \in S(\R),
$$
and the continuity of $p_{\lambda_0}(x,D)$ from  $H^{\psi_0,1}(\R)$ into $H^{\psi_0,(-1-2\sigma)}(\R)$ that for $k\rightarrow \infty$
$$
<p_{\lambda_0}(x,D)u,v> = B_{\lambda}(u,v)=<f,v>
$$
for all $v \in S(\R)$. Thus $p_{\lambda_0}(x,D)u = f$. The uniqueness follows of course once again from (\ref{D}).
\n In order to get more regularity for variational solutions or equivalently for solutions to (\ref{E}) we assume that for $\lambda \geq \lambda_0$ the function $p_{\lambda}^{-1}(x,\xi):= \frac{1}{p(x,\xi)+\lambda}$ belongs to $S^{-2+\tau_0,\psi_0}_{\rho}(\R)$ for some $\tau_0 >0$. In this case we can prove
\begin{theorem}
Let $p(x,\xi)$ be given by (\ref{N}) where we assume for $q$ condition (\ref{P}) and for $f$ we require (\ref{LL}), (\ref{M}) to hold. In addition we suppose that $p \in S^{2+\tau_1,\psi_1}_{\rho}(\R) \subset S^{2+\tau_1+2\sigma+\tau_1\sigma,\psi_0}_{\rho}(\R)$ and $p^{-1}_{\lambda} \in S^{-2+\tau_0,\psi_0}_{\rho}(\R)$, \n$\tau_1+\tau_0 +2\sigma+ \tau_1\sigma<1$. Let $ u \in H^{p_{\lambda_0}}(\R) \subset H^{\psi_0,1}(\R)$ be the solution to (\ref{E}) for $f \in H^{\psi_0,k}(\R)$, $k \geq 0$. Then it follows that $u \in H^{\psi_0,2+k-\tau_0}(\R)$.
\end{theorem}
\begin{proof}
From Theorem 1.8 it follows that 
\begin{equation}
p^{-1}_{\lambda_0}(x,D) \circ p_{\lambda_0}(x,D)= id +r(x,D)
\label{Q}
\end{equation}
with $r \in S_0^{-1+\tau_1+\tau_0+2\sigma+ \tau_1\sigma,\psi_0}(\R)$. Since $p_{\lambda_0}(x,D)u=f$ we deduce from (\ref{Q}) that
\begin{eqnarray*}
u & = & p^{-1}_{\lambda_0}(x,D) \circ p_{\lambda_0}(x,D)u -r(x,D)u \\
& = & p^{-1}_{\lambda_0}(x,D)f - r(x,D)u.
\end{eqnarray*}
Now, $p^{-1}_{\lambda_0}(x,D)f \in H^{\psi_0,k+2-\tau_0}(\R)$ and $r(x,D)u \in H^{\psi_0,2-\tau_1-\tau_0-2\sigma-\tau_1\sigma}(\R)$ implying that $u \in H^{\psi_0,t}(\R)$ for $t =(k+2-\tau_0)\wedge(2-\tau_1-\tau_0-2\sigma-\tau_1\sigma)>1 $. With a finite number of iterations we arrive at $u \in H^{\psi_0,2+k-\tau_0}(\R)$.
\end{proof}
\begin{remark}
From $\tau_1+\tau_0+2\sigma+\tau_1\sigma<1$ the necessary condition $\sigma <\frac{1}{2}$ follows.
\end{remark}
\begin{corollary}
In the situation of Theorem 2.1, if $2+k-\tau_0 > \frac{n}{2\rho_0 \rho_1}$, compare (\ref{R}), then $u \in C_{\infty}(\R)$.
\end{corollary}
Finally we can collect all preparatory material to prove
\begin{theorem}
Let $f:\R\times [0,\infty) \rightarrow \mathbb{R}$ be an arbitrarily often differentiable function such that for $y \in \R$ fixed, the function $s\rightarrow f(y,s)$ is a Bernstein function. Moreover assume (\ref{LL}), (\ref{M}) and (\ref{GG}). In addition let $\psi:\R \rightarrow \mathbb{R}$ be a continuous negative definite function in the class $\Lambda$ which satisfies in addition (\ref{FF}). For an elliptic symbol $q \in S^{2,\psi}_{\rho}(\R)$ satisfying (\ref{P}) we define $p(x,\xi)$ by (\ref{N}). For $\psi_1$ and $\psi_2$ defined by (\ref{AAB}) and (\ref{AAC}), respectively we assume (\ref{JJ}). Suppose that $p \in S^{2 +\tau_1,\psi_1}_{\rho}(\R)$ and $\frac{1}{p+\lambda} \in S^{-2+\tau_0,\psi_0}_{\rho}(\R)$. If $\tau_1+\tau_0 +\sigma(2+\tau_1) < 1$, $\sigma$ as in (\ref{JJ}), then $-p(x,D)$ extends to a generator of a Feller semigroup on $C_{\infty}(\R)$.
\end{theorem}
\begin{proof}
We want to apply the Hille-Yosida-Ray theorem, compare [\ref{AAL}], Theorem 4.5.3. We know that $p(x,D)$ maps $H^{\psi_0,2+k+2\sigma+\tau_1+\tau_1 \sigma}(\R)$ into $H^{\psi_0,k}(\R)$. Hence if $k>\frac{n}{2\rho_0 \rho_1}$ the operator $(-p(x,D),H^{\psi_0,2+k+2\sigma+\tau_1+\tau_1\sigma}(\R))$ is densely defined on $C_{\infty}(\R)$ with range in $C_{\infty}(\R)$. That $-p(x,D)$ satisfies the positive maximum principle on $H^{\psi_0,2+k+2\sigma+\tau_1+\tau_1\sigma}(\R)$ follows from Theorem 2.6.1 in [\ref{AAM}]. Now, for $\lambda \geq \lambda_0$ we know that for $f \in H^{\psi_0,k+1}(\R)$ we have a unique solution to $p_{\lambda}(x,D)u=f$ belonging to $H^{\psi_0,2+k+1-\tau_0}(\R)$. But \n$\tau_1+\tau_0 +2\sigma+\tau_1\sigma < 1$ implies that $H^{\psi_0,2+k+1-\tau_0}(\R)\subset H^{\psi_0,2+k+2\sigma+\tau_1+\tau_1\sigma}(\R)$, hence for $f \in H^{\psi_0,k+1}(\R)$ we always have a (unique) solution \n$ u \in H^{\psi_0,2+k+2\sigma+\tau_1+\tau_1\sigma}(\R)$ implying the theorem.
\end{proof}
\section{Some Concrete Examples}
 The first part of this section will consider the work W.Hoh has done on \ps s with variable order of differentiation. We will consider the case where the Bernstein function $s\rightarrow f(s)$ is substituted by $(x,s)\rightarrow s^{r(x)}$ with $r:\R\rightarrow \mathbb{R}$ being a continuous function such that $0\leq r(x) \leq 1$ holds. Let $q:\R \times \R \rightarrow \mathbb{C}$ be a continuous function such that $\xi \rightarrow q(x,\xi)$ is a continuous negative definite function. It then follows that
 \begin{equation}
 \xi \rightarrow q(x,\xi)^{r(x)}
 \end{equation}
 is once again a continuous negative definite function implying that the \ps
 \begin{equation}
Au(x)  := -(2\pi)^{-\frac{n}{2}}\int_{\R}e^{ix\cdot \xi}q(x,\xi)^{r(x)}\hat{u}(\xi)d\xi
\end{equation}
is a candidate for a generator of a Feller semigroup. We now meet Hoh's result:
  
    \begin{theorem}
 Let $\psi:\R \rightarrow \mathbb{R}$ be a fixed continuous negative definite function such that its L\'{e}vy measure has a compact support and that
 \begin{equation}
 \psi(\xi)\geq c_0|\xi|^r, \ \ \ |\xi| \mbox{ large and } r>0,
 \end{equation}
 holds. Let $ q \in S^{2,\psi}_{\rho}(\R)$ be a real-valued negative definite symbol which is elliptic, i.e. we have
 \begin{equation}
 q(x,\xi)\geq \delta_0(1+\psi(\xi)).
 \end{equation}
 Further let $m:\R\rightarrow (0,1]$ be an element in $C^{\infty}_b(\R)$ satisfying
 \begin{equation}
 M-\mu <\frac{1}{2}
 \end{equation}
 where $M:= \sup m(x)$ and $0< \mu : = \inf m(x)$. Consider the symbol
 \begin{equation}
 (x,\xi) \rightarrow p(x,\xi):= q(x,\xi)^{m(x)}
 \end{equation}
 which has the property that $\xi \rightarrow p(x,\xi)$ is a continuous negative definite function. The operator
 \begin{equation}
 -p(x,D)u(x):= -(2\pi)^{-\frac{n}{2}}\int_{\R}e^{ix\cdot \xi}p(x,\xi)\hat{u}(\xi)d\xi
 \end{equation}
 maps $C^{\infty}_0(\R)$ into $C_{\infty}(\R)$, is closeable in $C_{\infty}(\R)$ and its closure is a generator of a Feller semigroup.
\end{theorem}
\vskip5pt
 For a proof see W.Hoh [\ref{AAG}], compare also [\ref{AAF}].
\n We are now going to consider a further example. First note that the function $s\rightarrow \sqrt{s}(1-e^{-4\sqrt{s}})$ is a Bernstein function. Hence, using Corollary 3.9.36 in [\ref{AAL}], it follows that for $0 \leq \alpha \leq 1$ the function $s \rightarrow s^{\frac{\alpha}{2}}(1-e^{-4s^{\frac{\alpha}{2}}})$ is also a Bernstein function. Thus, given a negative definite symbol $q \in S^{2,\psi}_{\rho}(\R)$ we may consider the new symbol 
$$p(x,\xi)=(1+q(x,\xi))^{\frac{\alpha(x)}{2}}(1-e^{-4(1+q(x,\xi))^{\frac{\alpha(x)}{2}}})
$$ 
for $\alpha(\cdot)$ being an appropriate function.

 \begin{lemma}
 Let $q \in S^{2,\psi}_{\rho}(\R)$ be a real-valued negative definite symbol which is elliptic, i.e.
 $$
 q(x,\xi) \geq \delta_0(1+\psi(\xi)).
 $$
 Also let $\alpha(\cdot) :\R \rightarrow(0,1]$ be an element in $C^{\infty}_b(\R)$ satisfying
 $$
 m-\mu <\frac{1}{2}
 $$
 where $m=\sup{\frac{\alpha(x)}{2}}$ and $ \mu = \inf {\frac{\alpha(x)}{2}} >0$.
 \n Now if we let $p(x,\xi)=(1+q(x,\xi))^{\frac{\alpha(x)}{2}}(1-e^{-4(1+q(x,\xi))^{\frac{\alpha(x)}{2}}})$, then we have for all $\epsilon >0$ the estimates
 \begin{equation}
 |\parxi \parx p(x,\xi)| \leq c_{\alpha,\beta,\epsilon}p(x,\xi)(1+\psi(\xi))^{\frac{-\rho(|\alpha|)+\epsilon}{2}}
 \label{AAA}
 \end{equation}
 i.e. $p \in S^{2m+\epsilon,\psi}_{\rho}(\R).$
 
 \end{lemma}
 \begin{proof}
 We have to estimate
 \begin{eqnarray*}
 \parxi \parx p(x,\xi) & = & \parxi \parx ((1+q(x,\xi))^{\frac{\alpha(x)}{2}}(1-e^{-4(1+q(x,\xi))^{\frac{\alpha(x)}{2}}})) \\
 & = & \parxi \parx (e^{\frac{\alpha(x)}{2}\log (1+q(x,\xi))}(1-e^{-4(1+q(x,\xi))^\frac{\alpha(x)}{2}})).
 \end{eqnarray*}
 Using (2.19) in [\ref{AAL}] we get
 $$
 \parxi \parx (e^{\frac{\alpha(x)}{2}\log (1+q(x,\xi))}(1-e^{-4(1+q(x,\xi))^\frac{\alpha(x)}{2}})) =  
 $$
 $$
  \sum_{\alpha' \leq \alpha} \sum_{\beta' \leq \beta} \left ( \begin{array}{c} \alpha \\ \alpha' \end{array} \right ) \left ( \begin{array}{c} \beta \\ \beta' \end{array} \right ) (\partial^{\alpha'}_{\xi} \partial^{\beta'}_{x} e^{\frac{\alpha(x)}{2}\log (1+q(x,\xi))}) 
  $$
  \begin{equation} \times (\partial^{\alpha -\alpha'}_{\xi} \partial^{\beta-\beta'}_x (1-e^{-4(1+q(x,\xi))^\frac{\alpha(x)}{2}})). 
  \label{F}
 \end{equation}
 First consider
 $$
 |(\partial^{\alpha'}_{\xi} \partial^{\beta'}_{x} e^{\frac{\alpha(x)}{2}\log (1+q(x,\xi))})|.
 $$
 By (2.28) in [\ref{AAL}] with $l=|\alpha'|+|\beta'|$ we get

 $$
  |(\partial^{\alpha'}_{\xi} \partial^{\beta'}_{x} e^{\frac{\alpha(x)}{2}\log (1+q(x,\xi))})| \leq 
  $$
   \begin{equation}
     e^{\frac{\alpha(x)}{2}\log(1+q(x,\xi))}   \sum_{\begin{array}{c} \alpha'^1 + \ldots + \alpha'^{l'}  = \alpha' \\ \beta'^1+ \ldots +\beta'^{l'}=\beta' \\ l'=0,1,\ldots , l \end{array}} |c_{\{\alpha'^j,\beta'^j\}} \prod_{j=1}^{l'} q_{\alpha'^j \beta'^j}(x,\xi) |,
  \label{G}
  \end{equation}
  where 
  \begin{eqnarray*}
  q_{\alpha'^j \beta'^j}(x,\xi) & = & \partial^{\alpha'^j}_{\xi} \partial^{\beta'^j}_x(\frac{\alpha(x)}{2} \log (1+q(x,\xi))) \\ & = & \sum_{\bar{\beta}'^j \leq \beta'^j} \left ( \begin{array}{c} \beta'^j \\ \bar{\beta}'^j \end{array} \right ) \left ( \partial^{\beta'^j - \bar{\beta}'^j}_x \frac{\alpha(x)}{2} \right ) \ \partial^{\alpha'^j}_{\xi} \partial^{\bar{\beta}'^j}_x \log(1+q(x,\xi)).
  \end{eqnarray*}
  Now, using (2.26) in [\ref{AAL}] with $k=|\alpha'^j|+|\bar{\beta}'^j|>0$ we get
  $$
  \partial^{\alpha'^j}_{\xi} \partial^{\bar{\beta}'^j}_x \log(1+q(x,\xi))=
  $$
  $$
   \sum_{\begin{array}{c} \tilde{\alpha}'^1+\ldots +\tilde{\alpha}'^{l'} \\ \tilde{\beta}'^1 + \ldots +\tilde{\beta}'^{l'}=\bar{\beta}'^j \end{array}} c_{\{\tilde{\alpha}'^j, \tilde{\beta}'^j \}} \prod_{i=1}^k \frac{\partial^{\tilde{\alpha}'^i}_{\xi} \partial^{\tilde{\beta}'^i}_x (1+q(x,\xi))}{(1+q(x,\xi))}.
  $$
  Since we assume that $q(x,\xi)$ is an elliptic symbol in $S^{2,\psi}_{\rho}(\R)$, we get 
  $$
  \left| \partial^{\alpha'^j}_{\xi} \partial^{\bar{\beta}'^j}_x \log(1+q(x,\xi)) \right|
   $$
  \begin{eqnarray*}
  & \leq &  c_{\alpha'^j, \bar{\beta}'^j} \sum_{\begin{array}{c} \tilde{\alpha}'^1+\ldots +\tilde{\alpha}'^{l'} \\ \tilde{\beta}'^1 + \ldots +\tilde{\beta}'^{l'}=\bar{\beta}'^j \end{array}} \prod_{i=1}^k (1+\psi(\xi))^{\frac{-\rho(|\tilde{\alpha}'^i|)}{2}}\\ & \leq  & c_{\alpha^j, \bar{\beta}^j}(1+\psi(\xi))^{\frac{-\rho(|\alpha'^j|)}{2}},
 \end{eqnarray*}
   where we used the subadditivity of $\rho$.We always have 
   $$|\log(1+q(x,\xi))| \leq c_{\epsilon}(1+\psi(\xi))^{\frac{\epsilon}{2l}}.$$
    It follows for $\alpha \in C^{\infty}_{b}(\R)$ that 
   \begin{equation}
   | q_{\alpha'^j,\beta'^j}(x,\xi)| \leq c_{\alpha'^j,\beta'^j,\epsilon} \left\{ \begin{array}{cc} 
   (1+\psi(\xi))^{\frac{-\rho(|\alpha'^j|)}{2}}, & \alpha'^j \neq 0 \\ (1+\psi(\xi))^{\frac{\epsilon}{2l}}, & \alpha^j=0\ . \end{array} \right. 
   \label{H}
   \end{equation}
   Putting (\ref{G}) and (\ref{H}) together we get 
   \begin{equation}
    |(\partial^{\alpha'}_{\xi} \partial^{\beta'}_{x} e^{\frac{\alpha(x)}{2}\log (1+q(x,\xi))})| \leq c_{\alpha',\beta',\epsilon}e^{\frac{\alpha(x)}{2}\log (1+q(x,\xi))}(1+\psi(\xi))^{\frac{-\rho(|\alpha'|)+\epsilon}{2}}.
    \label{L}
    \end{equation}
    For the desired result we need
    $$
     |\partial^{\alpha -\alpha'}_{\xi} \partial^{\beta-\beta'}_x (1-e^{-4(1+q(x,\xi))^\frac{\alpha(x)}{2}})| 
     $$
     $$
     \leq c_{\alpha', \beta', \alpha,\beta, \epsilon}(1-e^{-4(1+q(x,\xi))^{\frac{\alpha(x)}{2}}})(1+\psi(\xi))^{-\frac{\rho(|\alpha-\alpha'|)}{2}}.
     $$
      When $\alpha-\alpha'=0$ and $\beta-\beta'=0$ there is nothing to prove.
  \n Otherwise, by  (2.28) in [\ref{AAL}] with $l_2 = |\alpha-\alpha'|+|\beta-\beta'|$, we get
     
      $$
       |\partial^{\alpha -\alpha'}_{\xi} \partial^{\beta-\beta'}_x (1-e^{-4(1+q(x,\xi))^\frac{\alpha(x)}{2}})| \leq  
       $$
        \begin{equation}
               e^{-4(1+q(x,\xi))^\frac{\alpha(x)}{2}}  \\  | \sum c_{\{(\alpha-\alpha')^j,(\beta-\beta')^j\}} \prod_{j=1}^{l'_2} q_{(\alpha-\alpha')^j (\beta-\beta')^j}(x,\xi)|,
               \label{J} 
               \end{equation}
       where the sum is such that
       $$
       _{\begin{array}{c} (\alpha -\alpha')^1+\ldots +(\alpha-\alpha')^{l'_2}=(\alpha-\alpha') \\ (\beta-\beta')^1+\ldots +(\beta-\beta')^{l'_2}=(\beta-\beta') \\ l'_2= 0,1, \ldots , l_2 \end{array} },
       $$
       and where
       $$
       q_{(\alpha-\alpha')^j(\beta-\beta')^j}(x,\xi)=\partial^{(\alpha-\alpha')^j}_{\xi} \partial^{(\beta-\beta')^j}_x(4(1+q(x,\xi))^{\frac{\alpha(x)}{2}}).
       $$
       Since $q(x,\xi)$ is in the symbol class $S^{2,\psi}_{\rho}(\R)$ we have the estimate
       $$
       |q_{(\alpha-\alpha')^j (\beta-\beta')^j}(x,\xi)| \leq \tilde{L}(1+q(x,\xi)) \mbox{ \  for all } (\alpha -\alpha')^j ,(\beta-\beta')^j \in \R,
       $$
       where $\tilde{L}(\lambda)$ is a suitable polynomial $\geq0$ which might depend on $(\alpha-\alpha')^j$ and $(\beta-\beta')^j$. 
     Now returning to (\ref{J}) we get
     $$
       |\partial^{(\alpha-\alpha')}_{\xi} \partial^{(\beta-\beta')}_x(1-e^{-4(1+q(x,\xi))^\frac{\alpha(x)}{2}} )| \leq \tilde{L}(1+q(x, \xi)) e^{-4(1+q(x,\xi))^\frac{\alpha(x)}{2}}       $$
       $$
         =\frac{4(1+q(x,\xi))^\frac{\alpha(x)}{2}}{1+4(1+q(x,\xi))^\frac{\alpha(x)}{2}} \ \cdot \frac{1+4(1+q(x,\xi))^\frac{\alpha(x)}{2}} {4(1+q(x,\xi))^\frac{\alpha(x)}{2}} \tilde{L}(1+q(x,\xi)) e^{-4(1+q(x,\xi))^\frac{\alpha(x)}{2}}
         $$
         $$
         \times (1+\psi(\xi))^{-\frac{\rho(|\alpha-\alpha'|)}{2}}(1+\psi(\xi))^{\frac{\rho(|\alpha-\alpha'|)}{2}}$$
         $$\leq \frac{4(1+q(x,\xi))^\frac{\alpha(x)}{2}}{1+4(1+q(x,\xi))^\frac{\alpha(x)}{2}}(1+\psi(\xi))^{-\frac{\rho(|\alpha-\alpha'|)}{2}}\cdot c_0  
         $$
         since
         $$
         \left|        \frac{1+4(1+q(x,\xi))^\frac{\alpha(x)}{2}} {4(1+q(x,\xi))^\frac{\alpha(x)}{2}}(1+\psi(\xi))^{\frac{\rho(|\alpha-\alpha'|)}{2}}\tilde{L}(1+q(x,\xi)) e^{-4(1+q(x,\xi))^\frac{\alpha(x)}{2}} \right| \leq c_0.
         $$     
         Now using (2.7) in [\ref{AAL}] i.e for all  $a\geq0$ and $t\geq0$ the estimate
         $$
         \frac{at}{1+at} \leq 1-e^{-at}  ,
         $$ we get
 \begin{equation}
         |\partial^{(\alpha-\alpha')}_{\xi} \partial^{(\beta-\beta')}_x(1-e^{-4(1+q(x,\xi))^\frac{\alpha(x)}{2}} )| \leq c_0(1-e^{-4(1+q(x,\xi))^\frac{\alpha(x)}{2}} )   (1+\psi(\xi))^{-\frac{\rho(|\alpha-\alpha'|)}{2}}
                  \label{K} 
         \end{equation}
         
       Substituting (\ref{L}) and (\ref{K}) into (\ref{F})
       $$
       | \parxi \parx (e^{\frac{\alpha(x)}{2}\log (1+q(x,\xi))}(1-e^{-4(1+q(x,\xi))^\frac{\alpha(x)}{2}}))|
       $$
       \begin{eqnarray*}
               & \leq &   \sum_{\alpha' \leq \alpha} \sum_{\beta' \leq \beta} \left ( \begin{array}{c} \alpha \\ \alpha' \end{array} \right )       \left ( \begin{array}{c} \beta \\ \beta' \end{array} \right )c_{\alpha',\beta',\epsilon}e^{\frac{\alpha(x)}{2}\log (1+q(x,\xi))}\\ &  & \times (1+\psi(\xi))^{\frac{-\rho(|\alpha'|)+\epsilon}{2}}(1-e^{-4(1+q(x,\xi))^\frac{\alpha(x)}{2}})  (1+\psi(\xi))^{-\frac{\rho(|\alpha-\alpha'|)}{2}}   \\ & \leq & c_{\alpha,\beta,\epsilon} e^{\frac{\alpha(x)}{2}\log (1+q(x,\xi))}(1-e^{-4(1+q(x,\xi))^\frac{\alpha(x)}{2}}) \\ &  & \times (1+\psi(\xi))^{\frac{-\rho(|\alpha|)+\epsilon}{2}}  \\ & \leq & c_{\alpha,\beta,\epsilon} p(x,\xi)(1+\psi(\xi))^{\frac{-\rho(|\alpha|)+\epsilon}{2}}.
       \end{eqnarray*}
       The proof now follows from the estimate $p(x,\xi) \leq (1+\psi(\xi))^m$.
       \n

           \end{proof}
\begin{lemma}
The function $p_{\lambda}^{-1}(x,\xi) = \frac{1}{p(x,\xi)+\lambda}$ belongs to the class $S^{-2\mu+\epsilon,\psi}_{\rho}(\R)$.
\end{lemma}
\begin{proof}
Using (2.27) in [\ref{AAL}]  we find with $l = |\alpha| +|\beta|$ that
$$
|\parxi \parx p^{-1}_{\lambda}(x,\xi)| \leq \frac{1}{p_{\lambda}(x,\xi)} 
\sum_{\begin{array}{c} \alpha^1 + \cdots +\alpha^l  =  \alpha \\ \beta^1 + \cdots +\beta^l  =  \beta \end{array}}c_{\{\alpha^j, \beta^j\}} \prod_{j=1}^l \left| \frac{\partial^{\alpha^j}_{\xi} \partial^{\beta^j}_x p_{\lambda}(x,\xi)}{p_{\lambda}(x,\xi)} \right|.
$$ 
For any $\epsilon > 0$ we find using (\ref{AAA})
$$
\left| \frac{\partial^{\alpha^j}_{\xi} \partial^{\beta^j}_x p_{\lambda}(x,\xi)}{p_{\lambda}(x,\xi)} \right| \leq \tilde{c}_{\alpha^j,\beta^j}(1+\psi(\xi))^{\frac{-\rho(|\alpha^j |)+\epsilon}{2}}
$$
and the ellipticity assumption of $p(x,\xi)$ together with the subadditivity of $\rho$ yields
$$
\left| \parxi \parx p^{-1}_{\lambda}(x,\xi)\right| \leq \tilde{c}_{\alpha,\beta,\epsilon}(1+\psi(\xi))^{-\mu}(1+\psi(\xi))^{\frac{-\rho(|\alpha|) +\epsilon}{2}}
$$
which proves the lemma.
\end{proof}

\vskip15pt
{\large{\bf{ADDRESS:}}}
\bigskip
\newline KRISTIAN P  EVANS, NIELS JACOB
\newline Department of Mathematics
\newline University of Wales Swansea
\newline Singleton Park
\newline Swansea SA2 8PP
\newline United Kingdom
\newline e-mail: N.Jacob@swan.ac.uk
\newline e-mail: makpe@swan.ac.uk

\end{document}